\documentclass[10pt,twoside,final]{amsart}

\usepackage[english]{babel}
\usepackage{graphicx,epstopdf,epsfig}
\usepackage{amsfonts,epsfig,fancyhdr,graphics,amsmath,amssymb}

\title{Products of strictly hyperbolic conjugacy classes  in symplectic groups}

\newtheorem{definition}{Definition}[section]
\newtheorem{theorem}{Theorem}[section]
\newtheorem{proposition}[theorem]{Proposition}
\newtheorem{lemma}[theorem]{Lemma}
\newtheorem{corollary}[theorem]{Corollary}
\newtheorem{remark}[theorem]{Remark}
\newtheorem{example}[theorem]{Example}

\newcommand {\ch }{\operatorname{char}}
\newcommand {\SpG }{\operatorname{Sp}}
\newcommand {\PSp }{\operatorname{PSp}}
\newcommand {\GL }{\mathrm{GL}}

\newcommand {\PSL }{\mathrm{PSL}}
\newcommand {\PSU }{\mathrm{PSU}}

\newcommand {\Bahn }{\mathrm{B}}
\newcommand {\rank }{\operatorname{rank}}
\newcommand {\Fix }{\mathrm{Fix}}

\newcommand {\Neg }{\mathrm{Neg}}
\newcommand {\GF }{\mathrm{GF}}

\newcommand {\Idm }{\mathrm{I}}

\newcommand {\RevIdm }{\mathrm{R}}
\newcommand {\PC }{\mathrm{PC}}

\newcommand {\tr }{\operatorname{tr}}

\newcommand {\mr }{\mathrm{mr}}
\newcommand {\M }{\mathrm{M}}

\usepackage{ifdraft}
\usepackage{listlbls}
\ifdraft{\usepackage[outer]{showlabels}
	\usepackage{todonotes}}{}
	\ifdraft{\usepackage[outer]{showlabels}
		\usepackage{todonotes} \usepackage{datetime}}{}

\usepackage[pagebackref,colorlinks,citecolor=red,urlcolor=blue,linkcolor=green, bookmarks=false,hypertexnames=true]{hyperref}
\usepackage{orcidlink}
\usepackage{fancyhdr}
\usepackage{verbatim}

\begin{document}

\bibliographystyle{plain}

\setcounter{page}{1}

\thispagestyle{empty}


\author{Klaus Nielsen}\,\orcidlink{0009-0002-7676-2944}
\email{klaus@nielsen-kiel.de}

\ifdraft{\today \ \currenttime}{\date{November 15, 2024}}
\pagestyle{fancy}
\fancyhf{}
\fancyhead[OC]{Klaus Nielsen}
\fancyhead[EC]{Strictly hyperbolic conjugacy classes}
\fancyhead[OR]{\thepage}
\fancyhead[EL]{\thepage}

\maketitle

\begin{abstract}
We call a conjugacy class of the  symplectic group $\SpG(2n, K)$ over a field $K$ strictly hyperbolic if its minimal polynomial is of the form $q(x) q^*(x)$, where the polynomial $q(x)$ is prime to its reciprocal $q^*(x) := x^n q(x^{-1})$.
It is shown that the product of 2 cyclic, strictly hyperbolic conjugacy classes of $\SpG(2n, K)$ contains all nonscalar elements of $\SpG(2n, K)$.
It follows that the projective symplectic group has a conjugacy class of covering number 2, i.e. $\PSp(2n,K) = \Omega^2$ for some conjugacy class $\Omega$ of $\PSp(2n,K)$. This verifies a conjecture of J. G. Thompson in the special case
of a (finite) projective symplectic group.
\end{abstract}


\section{Introduction} \label{intro-sec}

 There is a conjecture, usally attributed to J. G. Thompson, that a finite simple nonabelian group $G$ has a conjugacy class $\Omega$ of covering number 2; i.e $G = \Omega^2$. 
 Thompson's conjecture has been verified for several groups, e.g. for the alternating  and the sporadic groups. 
 For references.see \cite{Lev-1999} or \cite{EllersGordeev-1998}.  
Furthermore, it is known that the following simple groups have a conjugacy class of covering number 2:
\begin{enumerate}
    \item $\PSL(n, K)$  (Lev \cite[Theorem 5]{Lev-1994} for $|K| \ge 4$,
     and  \cite[Theorem 3]{Lev-1999}) for all fields $K$, 
    \item $\PSp(2n, K)$ if $-1 \in K^2$,
    \item $\PSU(V, K, h)$ if $|K| > 9$ and $h$ is hyperbolic (Bünger \cite[Satz 4.1.19]{Bunger-1997}).
\end{enumerate}
 
 A short proof for Lev's results in the case $K \ne \GF(2)$ or $n$ odd was given in \cite{BN-1999}.
 That $\PSp(2n, K)$ has a conjugacy class of covering number 2 if $-1 \in K^2$, is a trivial consequence of a theorem of Wonenburger. She showed that  a symplectic transformation is the product of 2 skew-symplectic involutions \cite[Theorem 2]{Wonenburger-1966}. Hence if $-1 \in K^2$, then a symplectic transformation  the product of 2 symplectic skew-involutions.
 
A great breakthrough was made in 1998 by  
Ellers and  Gordeev \cite{EllersGordeev-1998}, who proved the Thompson conjecture for all simple groups of Lie type over a field with more than 8 elements.

Almost 30 years ago, I proved  that the projective symplectic group $\PSp(2n, K)$ has a conjugacy class of covering number 2 provided that $2n \ge 4$. But I lost interest in the subject, after Ellers and Gordeev had published their result.

Recently, Larsen and Tiep \cite{LarsenTiep-2025} showed that Thompson's conjecture holds for almost all finite simple groups. So there seems to be still some interest in Thompson's conjecture.

\begin{theorem} \label{theorem-1}
	Let $2n \ge 4$. Let $\Omega_1, \Omega_2$ be cyclic, strictly hyperbolic conjugacy classes of $\SpG(2n, K)$. 
	\begin{enumerate}
		\item If $K = \GF(3)$ or $2n \ge 6$ or $\Omega_1$ and $\Omega_2$ are triangularizable, then $\Omega_1 \Omega_2$  contains all nonscalar elements of $\SpG(2n, K)$.
		\item If $2n = 4$, then $\Omega_1 \Omega_2$  contains all elements $M$ of $\SpG(2n, K)$ with $M^2 \ne \Idm_4$.
		\end{enumerate}
\end{theorem}

Using the computer algebra system gap, we can show

\begin{example}              
	Let $\Omega, \Psi$ be cyclic conjuagcy classes of $\SpG(4,5)$ with $\mu(\Omega) = (x^2+2)(x^2 - 2), \mu(\Psi) =  (x-2)^2(x + 2)^2$. Then $\Omega \Psi$ contains no involution.
\end{example}

\begin{corollary} \label{cor-1}
	Let $2n \ge 4$, and let $2n \ge 6$ if $K = \GF(2)$.
	 Then $\SpG(2n, K) = \Omega^2 \cup -\Idm_{2n}$ for some cyclic, strictly hyperbolic conjugacy class of $\SpG(2n, K)$.
	 \end{corollary}

\begin{proof}
	By theorem \ref{theorem-1}, it suffices to show that there exists a monic polynomial $q \in K[x]$ of degree $n$ such that $q$ is prime to its reciprocal $q^*$.
	If $|K| \ge 4$, let $q(x) = (x-\lambda)^n$, where $\lambda \ne 0, \pm 1$.
	If $K = \GF(3)$ and $n=2m$ is even let $q(x) = (x^2-x-1)^m$.
	If $K = \GF(3)$ and $n=2m+3$ is odd let $q(x) = (x^2-x-1)^m (x^3 - x -1)$.
	If $K = \GF(2)$ let $q(x) = x^n + x^{m+1} + 1$, where $m = \lfloor \frac{n}{2} \rfloor$.
\end{proof}

For the proof of  \ref{theorem-1}, we prove 2  auxilary results.
Let $\SpG(2n, K)$  denote the standard symplectic group, the set of all automorphs of the matrix 
\[
G = \left (\begin{array} {cc} 0 & \Idm_n\\ - \Idm_n & 0 \end{array} \right );
\]
i.e. $P \in \SpG(2n, K)$ if $PGP' = G$, where $P'$ denotes the transpose of $P$.

Let $\Omega$ be conjugacy of $\SpG(2n, K)$. By $\PC(\Omega)$, we denote the set of all matrices $A \in \M(n, K)$ occuring as the upper principle corner of some 
\[
P = \left (\begin{array} {cc} A & B\\ C & D \end{array} \right )\in \Omega.
\]
It is an interesting open problem to determine $\PC(\Omega)$. In the general linear group the analogue  problem is completely solved by the interlacing theorem of R.C. Thompson \cite[Theorem 6]{Thompson-1979} and E. M. de S\'a \cite[Theorem 5.4]{deSa-1979}:

\begin{proposition} \label{INTERLACING}
	Let $A \in \M(p, K), P \in \M(p+q, K)$ with invariant factors $i_k(A)$ and $i_k(P)$. Then there exist matrices $X, Y, Z$ such that
	\[
	P \sim  \left (\begin{array} {cc} A & X\\ Y & Z \end{array} \right )
	\]
	if and only if $i_k(P) | i_k(A) |i_{k+2q}(P)$.
\end{proposition}

In the case of a conjugacy class of a hyperbolic unitary vector space, Bünger has shown that $\PC(\Omega)$ contains nonscalar matrix with determinant $\alpha$ provided that $\det \Omega =  \alpha \overline{\alpha}^{-1}$.

For our purposes, it is sufficient to show

\begin{lemma}   \label{lemma1}
	Let  $2n \ge 4$, and let $\Omega$ be a nonscalar conjugacy of $\SpG(2n, K)$. 
	\begin{enumerate} 
		\item Let $\ch K = 2$. Then $\PC(\Omega)$ contains
    	\begin{enumerate}
		\item the identity matrix $\Idm_n$, 
		\item a transvection or a  diagonal matrix $\Idm_{n-2} \oplus [\epsilon, \epsilon^{-1}]$, where $\epsilon \ne 1$,
		\item a diagonal matrix $\delta \Idm_1 \oplus \Idm_{n-1}$ or a matrix $Q \oplus \Idm_{n-2}$, where $Q$ has minimal polynomial $x^2+\delta$.
	\end{enumerate}
\item Let $\ch K \ne 2$, and assume that $i_2(\Omega) = x-1$. Then $\PC(\Omega)$ contains
		 a transvection and all diagonal matrices $\Idm_{n-1} \oplus \delta \Idm_1$.  
		\item Let $\ch K \ne 2$ and $i_2(\Omega) = 1$. Then $\PC(\Omega)$ contains
		 a nonprimary matrix with prescribed determinant $\delta$ if 
		 \begin{enumerate}
		 \item $2n \ge 6$ or 
		 \item $\Omega$ is not involutory and $K \ne \GF(3)$ or
		  \item $\Omega$ is not involutory and $\delta \ne 1$.
		\end{enumerate}
\item Let $\ch K \ne 2$, $2n = 4$, and let $\Omega \ne \pm \Idm_4$ be involutory.
      Then $\PC(\Omega)$ contains a skew-involution $A$ ($A^2 = -\Idm_2$).
	\end{enumerate}
\end{lemma}

We also need 
\begin{lemma} \label{lemma2}
	Let $n \ge 2$.
	Let $\Phi, \Delta$ be cyclic conjugacy classes of $\GL(n, K)$. Then $\Phi \Delta$ contains
	\begin{enumerate}
		\item all nonprimary matrices $M \in \GL(n, K)$ with determinant $\det \Phi \Delta$,
		\item a transvection if $\det \Phi \Delta$ = 1 and $\Delta \ne \Phi^{-1}$.
		\end{enumerate}
\end{lemma}
This is a result of Bünger and the author \cite[Theorem 4.2]{BN-1999} and generalizes an earlier result of Lev \cite[Lemma 5]{Lev-1993}.
For the convenience of the reader, we provide a short proof.

Using \ref{lemma1} and \ref{lemma2}, we can prove theorem \ref{theorem-1}:

\begin{proof}[Proof of theorem \ref{theorem-1}]
Let $\Omega$ be a nonscalar conjugacy class of $\SpG(2n, K)$. Let $\mu(\Omega_1) = q_1 q_1^*, \mu(\Omega_2) = q_2 q_2^*$, where $q_1$ and $q_2$ are prime to their reciprocals. Then 
$\PC(\Omega)$ contains a matrix $A = A_1 A_2$, where $A_1$ and $A_2$ cyclic, and $\mu(A_1) = q_1, \mu(A_2) = q_2$. 
If $K = \GF(3)$, $2n = 4$, and $\Omega$ is involutory, then 
 $\PC(\Omega)$ contains a skew-involution $A$, and 
\[
A = \left (\begin{array} {cc} 0 & 1\\ -1 & 0 \end{array} \right )
= \left (\begin{array} {cc} 1 & -1\\ 1 & 1 \end{array} \right )
\left (\begin{array} {cc} 1 & -1\\ 1 & 1 \end{array} \right )
\]
is the product the product of 2 matrices with minimal polynomial $x^2 + x - 1$.
Now let 
\[
W = \left (\begin{array} {cc} A & B\\ C & D \end{array} \right ).
\]
Then the Schur complement $[W/A] := D - CA^{-1} B$ of $A$ in $W$ is $A^+ = (A')^{-1} = A_1^+ A_2^+$. Hence
\[
W = \left (\begin{array} {cc} A & B\\ C & D \end{array} \right )
= \left (\begin{array} {cc} A_1 & 0\\ CA_2^{-1} & A_1^+ \end{array} \right )
\left (\begin{array} {cc} A_2 & A_1^{-1}B\\ 0 & A_1^+ \end{array} \right ) \in \Omega_1 \Omega_2.
\]
\end{proof}
In a similar way, one can deal with unitary and orthogonal groups.
Bünger  uses a theorem of Lev \cite[Theorem 2]{Lev-1994} to prove his result \cite[Satz 4.1.19]{Bunger-1997} in the unitary group $\PSU(V, K, h)$ of a hyperbolic  hermitian form $h$. Therefore he has  to assume that $|K| > 9$. To cover the remaining cases $K = \GF(4)$ and $K = \GF(9)$, one needs a new factorization for the general linear group. We will prove such a theorem in a following paper.

\section{Preliminaries}
 
 \subsection{Notations}
 \mbox{}

 For a linear mapping $\varphi$ of $V$ let $\Bahn^j(\varphi)$ denote the image and $\Fix^j(\varphi)$ the kernel of $(\varphi -1)^j$. The space $\Bahn(\varphi) := \Bahn^1(\varphi)$ is the path or residual space of $\varphi$.
 And $\Fix(\varphi) := \Fix^1(\varphi)$ is the fix space of $\varphi$. The space $\Neg(\varphi) = \Fix(-\varphi)$ is the negative space of  $\varphi$.
 
 By $\mu(\varphi)$ we denote the minimal polynomial of 
 $\varphi$.
 Clearly, the spaces $\Bahn^j(\varphi)$ and $\Fix^j(\varphi)$ are $\varphi$-invariant.

 \begin{remark} \label{remark-1}
 	Let $\varphi \in \SpG(V)$.
 	\begin{enumerate}
 		\item $\varphi$ is involutory if and only if $\Bahn(\varphi)  = \Neg(\varphi)$;
 		\item $\Bahn^j(\varphi)^{\perp} = \Fix^j(\varphi)$.
 	\end{enumerate}	
 \end{remark}
 
 \begin{definition}
 	 Let $\varphi$ be a linear transfromation of $V$. A subspace $T$ of $V$ is called anti-invariant (w.r.t. $\varphi)$ if $T \cap T \varphi = 0$. 
 	 The minimal rank $\mr(\varphi)$ of a linear transformation $\varphi$ is defined
 	 as $\mr(\varphi) := \min  \{\rank (\varphi - \lambda ); \lambda \in  \overline{K}\}$, where  $\overline{K}$ is an algebraic closure of $K$.
 \end{definition}
 
 We have $\mr(\varphi) = m$ if and only if $i_m(\varphi)= 1 \ne  i_{m+1}(\varphi)$. It is known that $\varphi$ has an anti-invariant subspace of dimension $m$ if and only 
 if $2m \le \dim V$ and $\mr(\varphi) \ge m$. This follows  e.g. immediately from
 \ref{INTERLACING}.
 
  \begin{lemma}                                                  \label{lemma2a}
 	Let $\varphi \in \SpG(V)$. If $\dim \Fix(\varphi) = n =: \frac{\dim V}{2}$, then $\varphi$ has a totally degenerate anti-invariant subspace of dimension $n$.
 \end{lemma}
 
 \begin{proof}
 It is easy to see that an arbitrary  subspace of dimension $\ge n$ has a totally 
 degenerate complement. So if $T$ is  a totally 
 degenerate complement of $\Fix(\varphi)$, then $T \cap \Bahn(\varphi)  = 0$. Hence 
 $V = T \oplus T(1-\varphi) = T \oplus T \varphi$.
 \end{proof}
 
 \subsection{Products of symmetric matrices}
 
 \begin{lemma}                                                  \label{lemma3}
 Let $P \in \M(n, K)$ be cyclic. Then there exists a symmetric matrix $S$ (with maximal Witt index) such that $P^S = P'$ is the transposed of $P$.
 \end{lemma}
 
 \begin{proof}
 	See e.g. \cite[§6]{Shoda-1929} or \cite[Theorem 66]{Kaplansky-1974}.
 \end{proof}

 \subsection{The symplectic matrix group}
 \mbox{}
 
 Usually, we work basis-free, but sometimes it is more convenient to compute with matrices. 
A matrix $P$ is an automorph of a matrix $G$ if $PGP' = G$.
By $\SpG(2n, K)$ we denote the standard symplectic group, the set of all automorphs of the skew-sum 
\[
G = [-\Idm_n \setminus \Idm_n] = \left (\begin{array} {cc} 0 & \Idm_n\\ -\Idm_n & 0 \end{array} \right ),
\]
  of $-\Idm_n$ and $\Idm_n$. Note that some authors and computer algebra systems like gap use the skew-sum of $-\RevIdm_n$ and $\RevIdm_n$ instead, where $\RevIdm_n$ denotes the  reverse identiy matrix of $\GL(n, K)$.

  \begin{lemma}                                              \label{lemma4}
  Let $\Omega$ be a conjugacy class of $\SpG(2n,K)$. Let 
  \[
  P  =  \left (\begin{array} {cc} A & B\\ C & D \end{array} \right ) \in \Omega.
  \]
  	Then
  	\begin{enumerate}
  		\item $A^X \in  \PC(\Omega)$ for all $X \in \GL(n, K)$.
  		\item $A + B S \in  \PC(\Omega)$ for all symmetric matrices $S \in \M(n, K)$.
  		\item $D \in  \PC(\Omega)$.
  	\end{enumerate}
  	
  	\begin{proof}
  		Let
  	\[
  	Y  = \left (\begin{array} {cc} X & 0\\ 0 & X^+ \end{array} \right ),
  	G  = \left (\begin{array} {cc} 0 & \Idm_n\\ -\Idm_n & 0 \end{array} \right ),
  	T  = \left (\begin{array} {cc} \Idm_n & 0\\ S & \Idm_n \end{array} \right ).
  	\]	
  	Then	
  	\[
  	P^X  =  \left (\begin{array} {cc} A^X & X^{-1}BX^+\\
  		X'CX & D^{X^+} \end{array} \right ),
  	P^G  =  \left (\begin{array} {cc} D & -C\\ -B & A \end{array} \right ),
  	\]
  	\[
  	P^T  =  \left (\begin{array} {cc} A+BS & B\\ C+DS -SA -SBS & D - SB \end{array} \right )  \in \Omega.
  	\]
  	\end{proof}
  \end{lemma}

 \subsection{Orthogonally indecomposable transformations}
 \mbox{}
 The proof of lemma \ref{lemma1} demands a closer look at orthogonally indecomposable
  symplectic transformations.
  
According to Huppert \cite[1.7 Satz]{Huppert-1980a}, an orthogonally indecomposable transformation of $\SpG(V)$ is either\footnote{Our type enumeration follows Huppert \cite{Huppert-1980a}. In \cite{Huppert-1990} Huppert uses  a different  enumeration.}
\begin{enumerate}
	\item bicyclic with elementary divisors $e_1 = e_2 =(x \pm 1)^m$
	(type 1) or
	\item indecomposable as an element of $\GL(V)$ (type 2) or
	\item cyclic with minimal polynomial $(qq^*)^t$, where $q$ is irreducible and prime to its reciprocal $q^*$, where $q^*(x) = q(0)^{-1}x^{\partial q}q(x^{-1})$ ( type 3).
\end{enumerate}

\begin{definition}
	Let $\varphi \in \SpG(V)$. We say that 
	$\varphi$ is hyperbolic if $V$ is the direct sum of 2 totally degenerate $\varphi$-invariant subspaces.
\end{definition}

\begin{definition}
	Let $\varphi \in \SpG(V)$ orthogonally indecomposable of type 1. We say that 
	$\varphi$ is orthogonally indecomposable of type 
	\begin{enumerate}
		\item 1o if $\dim V \equiv 2 \mod 4$,
		\item  1e if $\dim V \equiv 0 \mod 4$.
\end{enumerate}
\end{definition}

In \cite[2.3 Satz]{Huppert-1980a} and \cite[2.4 Satz]{Huppert-1980a}, Huppert proved the following: 

\begin{proposition} \label{INDECOMP1}
	Let $\ch K \ne 2$. 
	Let $\varphi \in \SpG(V)$ be orthogonally indecomposable  of type 1. Then $\dim V \equiv 2 \mod 4$, and $\varphi$ is hyperbolic.
\end{proposition}

\begin{corollary} \label{INDECOMP2}
	Let $\ch K \ne 2$.
	Let $\varphi \in \SpG(V)$ be orthogonally indecomposable  of type 1 with minimal polynomial $(x+1)^{2m+1}$. Then 
	$\varphi$ has a symplectic square root with minimal polynomial $(x^2+1)^{2m+1}$.
\end{corollary}

\begin{lemma}                                               \label{INDECOMP3}
	Let $\ch K = 2$. Let $\varphi \in \SpG(V, f)$ be orthogonally indecomposable of type 1e. Then $\varphi$ is hyperbolic. 
\end{lemma}

\begin{proof}
	See Bünger \cite[Satz 1.5.5]{Bunger-1997}.
	We give a short proof:
	
	Clearly, $\varphi$ is hyperbolic if $\dim V = 4$.
	So let $\dim V \ge 8$. Obviously, $\varphi$ is orthogonally indecomposable on $\Bahn(\varphi)/\Fix(\varphi)$. By induction, $\Bahn(\varphi)/\Fix(\varphi) = [\langle u \rangle_{\varphi} + \Fix(\varphi)]/\Fix(\varphi) \oplus [\langle w \rangle_{\varphi} + \Fix(\varphi)]/\Fix(\varphi)$, where $\langle u \rangle_{\varphi}$ and $\langle w \rangle_{\varphi}$ are totally degenerate.
	
	Let $a \in \Fix^3(\varphi)$.
	Then $\langle a \rangle_{\varphi}$  is totally degenerate:
	
	If $\dim V \ge 12$, then $\Fix^3(\varphi)$ is totally degenerate, if $\dim V = 8$, then $\langle a \rangle_{\varphi}$ is totally degenerate as 
	$\varphi$ is orthogonally indecomposable.
	
	Let $u = x + x\varphi$. 
	We compute $f_{i,j} := f((a+x)(1+\varphi)^i, (a+x)(1+\varphi)^j) = f(x(1+\varphi)^i, x(1+\varphi)^j) +  f(x(1+\varphi)^i, a(1+\varphi)^j) +  f(a(1+\varphi)^i, x(1+\varphi)^j)$, as $\langle a \rangle_{\varphi}$ is totally degenerate. Clearly, $f_{i,j} = 0$ if $i+j \ge 3$. Further $f_{0,2} = 
	f_{1,1} = 0$. Finally,  $f_{0,1} = f(x, x\varphi) + f(x,a\varphi) + f(a,x\varphi) = f(x, x\varphi) + f(x(1+\varphi)^2, a\varphi)$. Hence 
	$\langle x+a \rangle_{\varphi}$ is totally degenerate for a suitable $a$, as
	$\Fix^3(\varphi) \not \subseteq [x(1+\varphi)^2]^{\perp}$.
	Similarly, there exist $z \in V- \Bahn(\varphi), c \in \Fix^3(\varphi)$ such that $w = z + z \varphi$ and $\langle z+c \rangle_{\varphi}$ is totally degenerate. Obviously, $V = \langle x+a \rangle_{\varphi} \oplus \langle z+c \rangle_{\varphi}$.
\end{proof}


\subsubsection{Type 2, 3}

\begin{lemma}                                                  \label{SKEWINV0}
	Let  $\varphi \in \SpG(V, f)$ be orthogonally indecomposable of type 2 or 3. Then in a suitable basis of $V$, 
	\[
	\varphi = \left (\begin{array} {cc} 0 & B  \\ -B^{-1} & D \end{array} \right ),
	f = \left (\begin{array} {cc} 0 & \Idm_m  \\ -\Idm_m & 0 \end{array} \right ),
	\]
	where  $D$ is indecomposable, and $B$ is symmetric.
\end{lemma}

\begin{proof}
	It is well known that a linear transformation is the product of 2 
	involutions if and only if it is  similar to its inverse. 
	This is a result usually attributed to Wonenburger \cite[Theorem 1]{Wonenburger-1966} and Dokovi\'c \cite[Theorem 1]{Dokovic-1967}. The same authors have also shown that a symplectic transformation is the product of 2 skew-symplectic involutions, Wonenburger \cite[Theorem 2]{Wonenburger-1966} for $\ch K \ne 2$ and Dokovi\'c \cite[Theorem 1]{Dokovic1971} for $\ch K = 2$.
	In fact, an involution $\rho \in \GL(V)$ inverting $\varphi$ must  already be skew-symplectic:
	Let $V = \langle u \rangle_{\varphi}$, and put $w = u \rho$. Then
	$f(u \varphi^i \rho , w \varphi^j \rho) = f(w \varphi^{-i}, u \varphi^{-j}) = f(w \varphi^j, u \varphi^i) = -f(u \varphi^i, w \varphi^j)$.
	
	So let $\varphi = \sigma \tau$ be the the product of 2 skew-symplectic involutions. If  $\ch K \ne 2$ or $\varphi$ is not unipotent, then $\Fix(\sigma)$ and $\Neg(\sigma)$ are totally degenerate, and $\dim \Fix(\sigma) = \dim \Neg(\sigma) = n$.
	
	Assume first that   $\Neg(\varphi) = 0$ or $\ch K \ne  2$. Replacing $\varphi$ by $-\varphi$, we may assume that $\Neg(\varphi) = 0$.
	Then $V = \Fix(\sigma) \oplus \Fix(\sigma) \tau$ or $V = \Neg(\sigma) \oplus \Neg(\sigma) \tau$. 
	
	So let $\ch K = 2$, and let $\varphi$ be unipotent. Then $\varphi = \sigma \tau$  is the product of 2 involutions with $\dim \Fix(\sigma) = n+1 = \dim \Fix(\tau) +1$: $\varphi$ is similar to 
	
	\[
	\varphi = \left (\begin{array} {cc} 0 & \Idm_m \\ \Idm_m & N \end{array} \right )
		 = \left (\begin{array} {cc} \Idm_m &  0 \\ N & \Idm_m \end{array} \right )
	 \left (\begin{array} {cc} 0 & \Idm_m  \\ \Idm_m & 0 \end{array} \right ),
	\]
	where $N$ is nilpotent and cyclic.
	 Let $S$ be a complement of $\Fix(\varphi)$  in $\Fix(\sigma)$ containing $\Bahn(\sigma)$. Then $S$ is totally degenerate,
	and $V = S \oplus S\tau$.
	
	Hence in a suitable basis,
	\[
	\sigma = \left (\begin{array} {cc} \pm \Idm_m & 0 \\ D & \mp \Idm_m \end{array} \right ),
	\tau = \left (\begin{array} {cc} 0 & B \\ B^{-1} & 0 \end{array} \right ), 
	f = \left (\begin{array} {cc} 0 & \Idm_m  \\ -\Idm_m & 0 \end{array} \right ),
	\]
	Now $B$ must be symmetric as $\tau$ is skewsymplectic. Finally, $\varphi$ is similar to 
	\[
	\left (\begin{array} {cc} B^{-1} & 0  \\ 0 & \Idm_m \end{array} \right )
	 \left (\begin{array} {cc} 0 & B  \\ -B^{-1} & D \end{array} \right )
      \left (\begin{array} {cc} B & 0  \\ 0 & \Idm_m \end{array} \right )
      = \left (\begin{array} {cc} 0 & \Idm_m  \\ -\Idm_m & D \end{array} \right ).
	\]
	Hence $D$ must be indecomposable.
\end{proof}

\subsubsection{Type 1e}
\mbox{}

Following Kaplansky \cite[p.80]{Kaplansky-1974}, we say that a symplectic transformation $\varphi$ is alternate if $f(v \varphi, v) = 0$ for all $v \in V$.
If $\ch K \ne 2$, then all symplectic involutions are alternate.

\begin{lemma}                                               \label{SKEWINV1}
	 Let $\varphi \in \SpG(V, f)$. Let $\varphi = \sigma \tau$ be the product of 2 alternate involutions $\sigma, \tau$. If $U$ is a $\varphi$-cyclic subspace of $V$, then $U \sigma, U\tau  \le U^{\perp}$.
\end{lemma}

\begin{proof}
	Let $U = \langle u \rangle_{\varphi}$. It suffices to show that $u \sigma \in U^{\perp}$. We have $f(u\sigma, u\varphi^{2j}) =
	f(u\sigma \varphi^{-j}, u\varphi^j)= f(u \varphi^j \sigma, u \varphi^j) = 0$
	and $f(u\sigma, u\varphi^{2j+1}) = f(u\sigma \varphi^{-j}, u\varphi^{j+1}) = f(u \varphi^j \sigma, u\varphi^{j+1}) = f(u \varphi^j \sigma, u\varphi^j \sigma \tau) = 0$.
\end{proof}
	
\begin{lemma}                                               \label{SKEWINV2}
	Let $\ch K = 2$. Let 
	\[
	H = \left (\begin{array} {cc} 0 & \Idm_t\\ \Idm_t & 0  \end{array} \right ),
	M = \left (\begin{array} {cc} N_t & 0\\ 0 & N_t'  \end{array} \right ),
	P = \left (\begin{array} {cc} 0 & H\\ H & M  \end{array} \right ),
	\]
	where $N_t$ is a  cyclic nilpotent matrix.
	Then $P \in \SpG(2n, K)$ is orthogonally indecomposable of type 1e.
\end{lemma}

\begin{proof}
	\[
	P = \left (\begin{array} {cc} 0 & H\\ H & M  \end{array} \right )
 = \left (\begin{array} {cc} 0 & H\\ H & 0 \end{array} \right )
	\left (\begin{array} {cc} \Idm_n & HM\\ 0 & \Idm_n \end{array} \right )
	\]
	 is the product of 2 alternate involutions. Further, $P$ is bicyclic with minimal polynomial $(x+1)^n$.
	Hence  $P$ must be orthogonally indecomposable of type 1: otherwise $\dim \Fix(T) = n$ by \ref{SKEWINV1}.

\end{proof}

\begin{corollary} \label{cor-2}
	Let $P \in \SpG(2n, K)$. Then $P$ is conjugate to a matrix 
	\[
	P_1 := \left (\begin{array} {cc} 0 & B  \\ -B^{-1} & D \end{array} \right ),
	\]
	where $B$ is symmetric,
	 if and only if $P$ has no elementary divisor $(x \pm 1)^{2t+1}$ of odd degree.
\end{corollary}

\begin{proof}
	It follows from \ref{SKEWINV0} and \ref{SKEWINV2} that $P$ is conjugate to $P_1$ if $P$ has no elementary divisor $(x \pm 1)^{2t+1}$ of odd degree. Conversely, $P_1$ is similar to a matrix 
	\[
	M_D := \left (\begin{array} {cc} 0 & \Idm_n  \\ - \Idm_n & D \end{array} \right ).
	\]
	If $D = X \oplus Y$, then $M_D$ is similar to $M_X \oplus M_Y$. So we may assume that $K$ is algebraicly closed and that $D$ is cyclic with minimum polynomial $(x-\lambda)^n$. But then $M_D$ is cyclic with minimum polynomial $(x^2- \lambda x + 1)^n$.
\end{proof}

\begin{remark}
	The Dickson transform of a polynomial $p(x) \in K[x]$ of degree $m$ is the polynomial 
	$p(x+x^{-1})x^n$. It can be shown that the invariant divisors of $P_1$ are the Dickson transforms of the invariant divisors of $D$. We have just used this fact in the proofs  of 
	\ref{cor-2} and \ref{SKEWINV0}.
\end{remark}

\subsubsection{Type 1o}
\mbox{}

\begin{lemma}                                               \label{SKEWINV2a}
	Let $\ch K \ne 2$. Let $\varphi \in \SpG(V)$ be orthogonally indecomposable
	of type 1 with minimal polynomial $(x+1)^{2m+1}$. Then
	$\varphi$ is conjugate to a matrix
	\[
	\left (\begin{array} {cc} -\Idm_n & B \\ C & -\Idm_n + CB \end{array} \right ) \in \SpG(2n, K),
	\]
	where $\rank B = \rank C = n-1$.
\end{lemma}

\begin{proof}
	By \ref{INDECOMP2}, $\varphi$ has a symplectic square root with minimal polynomial $(x^2+1)^{2m+1}$. Using \ref{cor-2}, we see that $\varphi$ is conjugate to a matrix
	\[
	P = \left (\begin{array} {cc} 0 & B \\ -B^{-1} & D \end{array} \right ) 
	\left (\begin{array} {cc} 0 & B \\ -B^{-1} & D \end{array} \right )
	= \left (\begin{array} {cc} -\Idm_n & BD \\ -DB^{-1} & D^2 -\Idm_n  \end{array} \right ).
	\]
	As above, we observe that $D$ must be cyclic and nilpotent.
\end{proof}

\begin{lemma}                                                    \label{SKEWINV3}
	Let $\ch K = 2$. Let $\varphi \in \GL(W)$ be cyclic with
	characteristic polynomial $(x+1)^{2m+1}$, where $m\ge 1$. Let
	$\varphi = \sigma \tau$ be a product of two involutions. Then
	\begin{enumerate}
		\item $\dim \Fix(\sigma) = m+1 = \dim \Fix(\tau),
		\dim \Bahn(\sigma) = m = \dim \Bahn(\tau) $.
		\item $\Bahn(\varphi) = \Bahn(\sigma) \oplus \Bahn(\tau)$,
		$\Fix(\varphi) = \Fix(\sigma) \cap \Fix(\tau)$.
		\item If $\Bahn(\sigma) \cap \Fix(\varphi) = 0$, then \\
		$V = \Fix(\tau) \oplus \Bahn(\sigma)$,
		$\Bahn(\sigma) \not \subseteq \Bahn^2(\varphi)$,
		$\Fix(\varphi) \subseteq \Bahn(\tau) \subseteq \Bahn^2(\varphi)$.
		\item Either $\Bahn(\varphi) = \Bahn(\tau) \oplus \Bahn(\tau)\varphi$
		or $\Bahn(\varphi) = \Bahn(\sigma) \oplus \Bahn(\sigma)\varphi$.
		\item Either $V = \Fix(\tau) \oplus \Bahn(\sigma)$ or $V = \Fix(\sigma) \oplus \Bahn(\tau)$.
	\end{enumerate}
\end{lemma}

\begin{proof}
	1 and 2 are clear.
	We prove 3: Suppose that
	$\Bahn(\sigma) \cap \Fix(\varphi) = 0$. Then
	$\Fix(\sigma) = \Bahn(\sigma) \oplus \Fix(\varphi)\subseteq \Bahn(\varphi)$ and
	$\Bahn^2(\varphi) \supseteq     \Fix(\sigma)(1+\varphi) = \Bahn(\tau)$.
	It follows that $\Bahn(\sigma) \not \subseteq \Bahn^2(\varphi)$.
	From $\varphi = \sigma \tau = \tau \sigma^{\tau}$ we deduce that
	$\Fix(\varphi) \subseteq \Bahn(\tau)$.
\end{proof}

\begin{lemma}                                               \label{SKEWINV4}
Let $\ch K = 2$.
Let $\varphi \in \SpG(V)$ be orthogonally indecomposable
of type 1, and let $\dim V = 2n \equiv 2 \mod 4$.
Let $\varphi = \sigma \tau$ be a product of two involutions $\sigma, \tau \in \SpG(V)$.
Let $V = U_1 \oplus U_2$, where $U_1, U_2$ are
$\langle \sigma, \tau \rangle$-invariant. Then

\begin{enumerate}
	\item $\Bahn(\varphi) = \Bahn(\sigma) \oplus \Bahn(\tau)$.
	\item $\Fix(\varphi) = [\Bahn(\sigma) \cap \Fix(\tau)] \oplus
	[\Bahn(\tau) \cap \Fix(\sigma)]$.
	\item 
	\begin{enumerate}
	\item $U_1 = \Fix(\sigma_1) \oplus \Bahn(\tau_1)$ and 
	$U_2 = \Bahn(\sigma_2) \oplus \Fix(\tau_2)$ or 
	\item $U_1 = \Bahn(\sigma_1) \oplus \Fix(\tau_1)$ and
	$U_2 = \Fix(\sigma_2) \oplus \Bahn(\tau_2)$.
\end{enumerate}
	\item If $\sigma$ and $\tau$ are alternate, then $U_1$ and
	$U_2$ are totally degenerate.
\end{enumerate}
\end{lemma}

\begin{proof}
	1: follows from \ref{SKEWINV4}(2). \\
	2: We have $\Fix(\varphi) \not \subseteq \Bahn(\sigma)$: Otherwise
	$\Fix(\sigma) \subseteq \Bahn(\varphi)$, and by 1),
	$\Fix(\sigma) = \Fix(\sigma) \cap  \Bahn(\varphi) = \Bahn(\sigma)
	\oplus [ \Bahn(\tau) \cap \Fix(\sigma) ] = \Bahn(\sigma)$,
	a contradiction. Similarily,  $\Fix(\varphi) \not\subseteq \Bahn(\tau)$.
	3: follows from 2.\\
	4: If $\sigma$ and $\tau$ are alternate, then $U_1$ and
	$U_2$ are totally degenerate by \ref{SKEWINV1}.\\
\end{proof}

\begin{corollary}                                               \label{SKEWINV5}
	Let $\ch K = 2$. Let $\varphi \in \SpG(V)$ be orthogonally indecomposable
	of type 1o, and let $\dim V = 2n \equiv 2 \mod 4$.
	 Then
	  $\varphi$ is conjugate to a matrix
	 \[
	 \left (\begin{array} {cc} \Idm_n & B \\ C & \Idm_n + CB \end{array} \right ) \in \SpG(2n, K),
	 \]
	 where $\rank B = \rank C = n-1$.
\end{corollary}

\begin{proof}
We have a decomposition $V = S \oplus T$, with
$\Bahn(\sigma) \subseteq S \subseteq \Fix(\sigma)$ and $\Bahn(\tau) \subseteq T \subseteq \Fix(\tau)$:
Take  $S = \Fix(\sigma_1) \oplus \Bahn(\sigma_2), T = \Bahn(\tau_1) \oplus \Fix(\tau_2)$ or
$S = \Bahn(\sigma_1) \oplus \Fix(\sigma_2),
T = \Fix(\tau_1) \oplus \Bahn(\tau_2)$.
Then in a suitable basis
\[
\sigma = \left (\begin{array} {cc} \Idm_n & 0 \\ C & \Idm_n  \end{array} \right ),
\tau = \left (\begin{array} {cc} \Idm_n & B \\ 0 & \Idm_n  \end{array} \right ),
\varphi = \left (\begin{array} {cc} \Idm_n & B \\ C & \Idm_n + CB \end{array} \right ).
\]
\end{proof}

\begin{remark}
	It can be shown  that $\varphi$ has a cyclic square root even if $\ch K = 2$; see
	\cite[Corollary 3.7]{KNielsen-2024}. 
	So \ref{SKEWINV5} follows also from \ref{cor-2}.
\end{remark}

\begin{corollary}                                               \label{SKEWINV6}
	 Let $\ch K = 2$. Let $P \in \SpG(2n, K)$. Then $P$ is conjugate to a matrix
	\[
	\left (\begin{array} {cc} \Idm_n & B \\ C & \Idm_n + CB \end{array} \right ).
	\]
\end{corollary}

\begin{proof}
If $P$ has no elementary divisor $(x\pm 1)^{2t+1}$ of odd degree, then using
\ref{cor-2}, we see that $P$ is conjugate to 
\[
\left (\begin{array} {cc} \Idm_n & 0\\ -B^{-1} & \Idm_n \end{array} \right )
\left (\begin{array} {cc} 0 & B\\ -B^{-1} & D  \end{array} \right )
\left (\begin{array} {cc} \Idm_n & 0\\ B^{-1} & \Idm_n \end{array} \right )
= \left (\begin{array} {cc} \Idm_n & B\\  (D-2\Idm_n)B^{-1}& D- \Idm_n \end{array} \right ).
\]
\end{proof}

\section{Principal corners of a symplectic conjugacy class: Proof of \ref{lemma1} }

\subsection{The case $\ch K = 2$}

\begin{lemma} \label{PCORN0}
	Let $\ch K = 2$. Let $\delta \in K^*$. Let $\Omega \subseteq \SpG(4, K)$ be a nonscalar conjugacy class. Then $\PC(\Omega)$ contains 
	\begin{enumerate}
		\item  a transvection or a nonprimary with determinant one,
		\item  a matrix with characteristic polynomial $(x+1)(x+\delta)$ or $(x^2+\delta)$.
	\end{enumerate}
\end{lemma}

\begin{proof}
	It follows from \ref{SKEWINV6} that $\Omega$ contains a matrix
	\[
	\left (\begin{array} {cc} \Idm_2 & B\\ C & \Idm_2 + CB\end{array} \right ).
	\]
	By \ref{lemma4}, $\PC(\Omega)$ contains all matrices $\Idm_2 + BS$, where $S$ is symmetric.
	If $\rank B = 1$ or $B$ has discriminant 1, then $\PC(\Omega)$ contains a transvection:
	\[
	\left (\begin{array} {cc} 1 & 0\\ 0 & 1 \end{array} \right ) +
	\left (\begin{array} {cc} \alpha & 0\\ 0 & 0 \end{array} \right )
	\left (\begin{array} {cc} 0 & 1\\ 1 &  0 \end{array} \right )
	= 	\left (\begin{array} {cc} 1   & \alpha \\  0  &  1 \end{array} \right ).
	\]
	\[
	\left (\begin{array} {cc} 1 & 0\\ 0 & 1 \end{array} \right ) +
	\left (\begin{array} {cc} 0 & 1\\ 1 & 0 \end{array} \right )
	\left (\begin{array} {cc} 1 & 0\\ 0 &  0\end{array} \right )
	= 	\left (\begin{array} {cc} 1   & 0 \\  1  &  1  \end{array} \right ).
	\]
	\[
	\left (\begin{array} {cc} 1 & 0\\ 0 & 1 \end{array} \right ) +
	\left (\begin{array} {cc} 1 & 0\\ 0 & 1 \end{array} \right )
	\left (\begin{array} {cc} 1 & 1\\ 1 &  1 \end{array} \right )
	= 	\left (\begin{array} {cc} 0   & 1 \\  1  &  0  \end{array} \right ).
	\]
	If $B$ is alternating, then $\PC(\Omega)$ contains all traceless matrices:
	\[
	\left (\begin{array} {cc} 1 & 0\\ 0 & 1 \end{array} \right ) +
	\left (\begin{array} {cc} 0 & 1\\ 1 & 0 \end{array} \right )
	\left (\begin{array} {cc} \lambda & \mu\\ \mu &  \nu\end{array} \right )
	= 	\left (\begin{array} {cc} 1+\mu   &  \nu \\  \lambda  &  1+\mu  \end{array} \right ).
	\]
	Clearly,  $\PC(\Omega)$ contains all diagonal matrices if $B$  is diagonal.
\end{proof}

\begin{corollary}   \label{PCORN1}
	Let $\ch K = 2$.                               
	Let $\Omega \subseteq \SpG(2n, K)$ be a nonscalar conjugacy class.
	Let $\delta \in K^*$. Then $\PC(\Omega)$ contains 
	\begin{enumerate}
		\item the identity matrix $\Idm_n$, 
		\item a transvection or all diagonal matrices $\Idm_{n-2} \oplus D_2$,
		\item a diagonal matrix $\delta \Idm_1 \oplus \Idm_{n-1}$ or a matrix $Q \oplus \Idm_{n-2}$, 
		where $Q$ has minimal polynomial $x^2+\delta$.
	\end{enumerate}	
\end{corollary}

\subsection{The case $\ch K \ne 2$}

\begin{lemma} \label{PCORN2}
	Let $\ch K \ne 2$. Let $\Omega \subseteq \SpG(4, K)$ be a nonscalar conjugacy class. Then
	$\PC(\Omega)$ contains 
	
	\begin{enumerate}
		\item a transvection if $i_2(\Omega) = x-1$, 
		\item  all traceless matrices of $\GL(2,K)$ if $\Omega$ is involutory,
		\item  all diagonal  matrices of $\GL(2,K)$ if $\Omega$ has no elementary divisor $x \pm 1$,
		\item  all nonscalar triangular  matrices of $\GL(2,K)$ if  $\Omega$ is not involutory and $\mu(\Omega) = (x - 1)p(x)$, where $p(1) \ne 0$. 
	\end{enumerate}
\end{lemma}

\begin{proof}
	We use \ref{lemma4} without further reference. If $\Omega$ consists of transvections, 
then $\Omega$ contains a matrix
\[
\left (\begin{array} {cc} \Idm_2 & B\\ 0 & \Idm_2 \end{array} \right ),
\]	
where $\rank B = 1$. Then $\PC(\Omega)$ contains 
\[
\left (\begin{array} {cc} 1 & 0\\ 0 & 1 \end{array} \right ) +
\left (\begin{array} {cc} \alpha & 0\\ 0 & 0 \end{array} \right )
\left (\begin{array} {cc} 0 & 1\\ 1 &  0 \end{array} \right )
= 	\left (\begin{array} {cc} 1   & \alpha \\  0  &  1  \end{array} \right ).
\]	
	
Now let $\mr(\Omega) \ge 2$.  By \ref{lemma2a} and \ref{cor-2}, $\Omega$ contains a matrix
\[
\left (\begin{array} {cc} 0 & B\\ -B^+ & D \end{array} \right ).
\]
2: If $\Omega$ is involutory, then $B$  is alternating, and $\PC(\Omega)$ contains all matrices $BS$, where $S \in \GL(2, K)$ is symmetric.\\
3: If   $\Omega$ has no elementary divisor $x \pm 1$, then we can assume that $B$ is symmetric. In this case, $\PC(\Omega)$ contains all diagonal matrices. \\
4: By \ref{cor-2}, $B$ be not symmetric. Furthermore,  $B \ne  -B'$: Otherwise, 
\[
\left (\begin{array} {cc} 0 & B\\ -B^+ & D \end{array} \right ) = 
\left (\begin{array} {cc} B & 0\\ 0 & B^{-1} \end{array} \right )
\left (\begin{array} {cc} 0 & \Idm_2\\ \Idm_2 & BD \end{array} \right )
\]
is the product of an alternating and a symmetric matrix. But then
$\Omega = - \Omega$, a contradiction. Hence we can assume that $B$ is (upper) triangular.
In this case, $\PC(\Omega)$ contains all nonsingular  (upper) triangular matrices.
\end{proof}

\begin{lemma}   \label{PCORN3}
	Let $\ch K \ne 2$, and $2n \ge 6$.                              
	Let $\Omega \subseteq \SpG(2n, K)$ be a nonscalar conjugacy class of $\SpG(2n, K)$.
	Let $\delta \in K^*$.
	Then $\PC(\Omega)$ contains 
	\begin{enumerate}
		\item  $\Idm_n$, a transvection, and all dilatations if $\Omega$ consist of transvections, 
		\item a nonprimary matrix with prescribed determinant if $\mr(\Omega) \ge 2$.
	\end{enumerate}	
\end{lemma}

\begin{proof}
\begin{enumerate}
	\item[a)] First assume that $\Omega$ has no elementary divisor $(x\pm1)^{2t+1}$ of odd degree.
	Then $\PC(\Omega)$ contains all diagonal matrices.
	\item[b)] Next assume that
	    $\Omega$ is orthogonally indecomposable of type 1
	    with $\mu(\Omega) = (x-\epsilon)^n$.
	    Then $\PC(\Omega)$ contains all diagonal matrices $D = \epsilon \Idm_1 \oplus D_2$
	   \item[c)] Let  $\Omega = \epsilon \Idm_2 \perp \Psi$, where $\Psi$ has no elementary divisor $(x\pm1)^{2t+1}$ of odd degree. Again, $\PC(\Omega)$ contains all diagonal matrices $D = \epsilon \Idm_1 \oplus D_2$.
	\end{enumerate}
	In any case, $\PC(\Omega)$ contains a nonprimary matrix with prescribed determinant.
 So we may assume that $2n = 6$ and $\Omega =  \Idm_2 \perp \epsilon \Idm_2 \perp \Psi$.
 \begin{enumerate}
    \item[i)] If $\Psi = \Idm_2$, then $\epsilon = -1$ and $\PC(\Omega)$ contains all matrices with mimimal polynomial $(x-1)(x^2+\delta)$.
    \item[ii)] If $\Psi = -\Idm_2$, then $\PC(\Omega)$ contains all matrices with mimimal polynomial $(x+\epsilon )(x^2+\delta)$.
    \item[ii)] Let $\Psi \ne \pm  \Idm_2$.
    \begin{enumerate}
     \item If $\mu_{\Psi}(1) \ne 0$, then $\PC(\Omega)$ contains all matrices  $\epsilon \Idm_1 \oplus U$, where $U$ is upper triangular.
     \item If $\epsilon = -1$ and $\mu_{\Psi}(-1) \ne 0$, then $\PC(\Omega)$ contains all matrices  $ \Idm_1 \oplus U$, where $U$ is upper triangular.
     \item If $\epsilon = 1$ and $\mu_{\Psi}(1) \ne 0$, then $\PC(\Omega)$ contains all matrices  $ \Idm_3$, all dialatations and transvections.
    \end{enumerate}
\end{enumerate}
\end{proof}

\section{Proof of lemma \ref{lemma2}}
\mbox{}

\begin{lemma}  \label{lemma5}
	Let $U \in \GL(m-1, K)$ is upper triangular. Then there exist polynomials $g_j(x)$ of degree $\le m-j -1$ with nonvanishing constant term such that the matrix 
	\[
	Y :=	\left (\begin{array} {cc}
		0 & U\\  y_0 & y  \end{array} \right ),
	\]
	where $y =(y_0, \dots, y_{m-1})$, has characteristic polynomial 
	\[ \chi(Y) = \det x \Idm_m - Y = x^m + \sum_{j=0}^{m-1} y_j g_j(x) x^j.
	\] 
	The polynomials $g_j(x) x^j$ are linearly independent.
\end{lemma}

\begin{lemma}                                                 \label{lemma6}
	Let $U_1 \in \GL(m,K), U_2 \in \GL(n-m,K)$ be upper 
	triangular matrices.  
	Let $\Phi$ be a cyclic similarity class of $\GL(n,K)$ with $\det \Phi = (-1)^{n+1} \det U_1 U_2$. 
	Then there exists a matrix $Z \in \M(m,n-m,K)$ such that 
	\[
	P_Z:= \left (\begin{array} {cc} 0 & \Idm_{n-1}\\ 1 & 0 \end{array} \right )
	\left (\begin{array} {cc} U_1 & Z\\ 0 & U_2 \end{array} \right )
	\in \Phi.
	\]
\end{lemma}

\begin{proof}
	Let $z = (z_{m-1}, \dots, z_1)', \widetilde{z} = (z_{m+1}, \dots. z_{n-1})$,
	\[
	U_1 = \left (\begin{array} {cc} u_1 & b\\ 0 & W_1 \end{array} \right ),
	U_2 = \left (\begin{array} {cc} u_2 & c\\ 0 & W_2 \end{array} \right ),
	Z = \left (\begin{array} {cc} z_m & y\\ z & 0 \end{array} \right ).	
	\]
	Then 
	\[
	P_Z = \left (\begin{array} {cc} 0 & \Idm_{n-1}\\ 1 & 0 \end{array} \right )
	\left (\begin{array} {cccc} u_1 & b & y_1 & y\\ 0 & W_1 & z & 0 \\ 
		0 & 0 & u_2 & c\\ 0 & 0 & 0 & W_2  \end{array} \right )
	=	\left (\begin{array} {cccc} 0 & W_1 & z & 0\\ 0 & 0 & u_2 & c \\ 
		0 & 0 & 0 & W_2\\ u_1 & b & z_m & \widetilde{z}  \end{array} \right ).
	\]
	Clearly, $P_Z$ is cyclic. Let 
	\[
	Y :=	\left (\begin{array} {cc}
		0 & W_2\\  z_m &  \widetilde{z} \end{array} \right ).
	\]
	We compute the characteristic polynomial $\chi(P_Z)$  of $P_Z$  by Laplace 
	expansion by the first m columns:
	\[
	\chi(P_Z) = x^m \chi(Y) + \sum_{j=1}^{m-1} z_j f_j(x)x^j,
	\]
	where the polynomials $f_j(x)$ are of degree $\le m-j-1$ and independent of $U_1$ and $U_2$, and $f_j(0) \ne 0$.
\end{proof}

\begin{lemma}                                                  \label{lemma7}
	Let $\Phi, \Delta$ be  similarity classes of cyclic 
	matrices in $\GL(n,K)$. Let $P_1 \in \GL(m,K), P_2 \in \GL(n-m,K)$ with $\det P_1 P_2 = \det \Phi \Delta$.
	Then there exists a matrix $C \in \M(m,n-m,K)$ such that
	\[
	\left (\begin{array} {cc} P_1 & C\\ 0 & P_2 \end{array} \right ) 
	\in \Delta \Phi.
	\]
\end{lemma}

  \begin{proof}  
  	We may assume that $P_1 = L_1 U_1$ and  $P_2 = L_2 U_2$, where $L_1, L_2$ are lower and  $U_1, U_2$ are upper triangular.    
  	By \ref{lemma6}, there exists a matrix $Z$ such that   
  	\[
  	\left (\begin{array} {cc} 0 & \Idm_{n-1}\\ 1 & 0 \end{array} \right )
  	\left (\begin{array} {cc} U_1 & Z\\ 0 & U_2 \end{array} \right )
  	\in \Delta.
  	\]
  	Further, there exists a matrix $X$ such that 
  	\[
  	\left (\begin{array} {cc} L_1 & X\\ 0 & L_2 \end{array} \right )
  	\left (\begin{array} {cc} 0 & 1\\ \Idm_{n-1} & 0 \end{array} \right )
  	\in \Omega:
  	\]
  	\[
  		\left (\begin{array} {cc} \RevIdm_m & 0\\ 0 &  \RevIdm_{n-m} \end{array} \right )
  	\left (\begin{array} {cc} L_1 & X\\ 0 & L_2 \end{array} \right )
  	\left (\begin{array} {cc} 0 & 1\\ \Idm_{n-1} & 0 \end{array} \right )
  		\left (\begin{array} {cc} \RevIdm_m & 0\\ 0 &  \RevIdm_{n-m} \end{array} \right )
  		\]
  		\[
  		\left (\begin{array} {cc} L_1^{\RevIdm_m} & \RevIdm_m X \RevIdm_{n-m}\\ 0 & L_2^{\RevIdm_{n-m}} \end{array} \right )
  		\left (\begin{array} {cc} 0 & 1\\ \Idm_{n-1} & 0 \end{array} \right ) \sim
  		\left (\begin{array} {cc} 0 & 1\\ \Idm_{n-1} & 0 \end{array} \right )
  		\left (\begin{array} {cc} L_1^{\RevIdm_m} & \RevIdm_m X \RevIdm_{n-m}\\ 0 & L_2^{\RevIdm_{n-m}} \end{array} \right ).
  		\]
  		Hence $\Phi \Delta$ contains the matrix
  		\[
  			\left (\begin{array} {cc} P_1 & L_1 X + Z U_2\\ 0 & P_2 \end{array} \right ) =
  		\left (\begin{array} {cc} L_1 & X\\ 0 & L_2 \end{array} \right )
  		\left (\begin{array} {cc} 0 & 1\\ \Idm_{n-1} & 0 \end{array} \right )
  		\left (\begin{array} {cc} 0 & \Idm_{n-1}\\ 1 & 0 \end{array} \right )
  		\left (\begin{array} {cc} U_1 & Z\\ 0 & U_2 \end{array} \right ) .
  	\]
\end{proof}

\begin{example}                                 \label{lemma8}
	Let $n \ge 2$. If $K = \GF(2)$ let $n \ge 3$.
	Let $\Phi, \Delta$ be cyclic similarity classes of  $\GL(n,K)$. 
	
	\begin{enumerate}
		\item $\Phi \Delta$ contains all matrices $M = X \oplus \delta \Idm_1$ if $\det M = \det \Phi \Delta$ and $\delta$ is not an eigenvalue of $X$.
		\item  If $\Phi$ is reducible, then $\Phi \Phi^{-1}$ contains  a transvection.
	\end{enumerate}
\end{example}

\begin{proof}
	1: Let $L \in \GL(n-1, K)$ be lower triangular, and let $U \in \GL(n-1, K)$ be upper triangular matrix. There exist matrices 
	\[
	A = \left (\begin{array} {cc} a & \delta \\ L & 0  \end{array} \right ) \in \Phi,
	B = \left (\begin{array} {cc} 0 & U\\ \epsilon & b  \end{array} \right ) \in \Delta.
	\]
	Then 
	\[
	\left (\begin{array} {cc} \delta \epsilon & aU + \delta b\\ 0 & LU  \end{array} \right ) \in \Phi \Delta.
	\]
	
2: There exist matrices $Q \in \GL(m, K), R \in \GL(n-m, K)$, and $N, M \in \M(m,n-m,K)$  such that  $\rank N,  \rank M, \rank(N - M) = 1$ and
\[
X = \left (\begin{array} {cc} Q & M\\ 0 & R  \end{array} \right ),
Y = \left (\begin{array} {cc} Q & N\\ 0 & R  \end{array} \right ) \in \Phi.
\]
 Now
\[
XY^{-1} = \left (\begin{array} {cc} \Idm_m & (M-N)R^{-1}\\ 0 & \Idm_{n-m}  \end{array} \right ) 
\]
is a transvection.
\end{proof}

\begin{remark}   
	It can be shown that $\Phi \Phi^{-1}$ contains  no transvection if $\Phi$ is irreducible. 
\end{remark}

\begin{lemma}               \label{ARADHERZOG}
	Let $\Phi, \Delta$ be nonscalar similarity classes of $\M(2,K)$. Assume that $\Phi$ is not irreducible.
	\begin{enumerate}
		\item  $\tr(\Phi \Delta) = K - \{\alpha \tr \Delta\}$ if and only if
		$(\Phi - \alpha \Idm_2)^2 = 0$ and $\Delta$ is irreducible.
		\item $\tr(\Phi \Delta) = K$ if and only if  $\Phi$ is not primary or $\Delta$ is not irreducible.
	\end{enumerate}
\end{lemma}

\begin{proof}
	Let  $\mu(\Phi) = (x - \lambda) (x - \alpha)$, $\mu(\Delta) = x^2 - \gamma x - \delta$. Then
	\[
	\left (\begin{array} {cc} \alpha & 0\\ \mu & \lambda \end{array} \right )
	\left (\begin{array} {cc} 0  & 1 \\ \gamma & \delta \end{array} \right ) =
	\left (\begin{array} {cc}  0 & \alpha   \\  \lambda \gamma   &  \lambda \delta + \mu \end{array} \right ) \in \Omega \Delta
	\]
	if $\lambda \ne \alpha$ or  $\mu \ne 0$.
	If $\lambda = \alpha$ and $\Psi$ is irreducible, then  $\tr(\Phi \Delta) = K - \{\alpha \tr \Delta\}$:
	\[
	\left (\begin{array} {cc} \alpha & \mu\\ 0 & \alpha \end{array} \right )
	\left (\begin{array} {cc} \delta - \epsilon & \rho \\ \nu & \epsilon \end{array} \right ) =
	\left (\begin{array} {cc} \alpha (\delta -\epsilon) + \mu \nu & \alpha \rho + \mu \epsilon \\  \alpha \nu   &  \alpha \epsilon  \end{array} \right ).
	\]
\end{proof}


\listoflabels

\end{document}

\typeout{get arXiv to do 4 passes: Label(s) may have changed. Rerun}